\documentclass[12pt]{amsart}
\usepackage{amsthm}
\newtheorem*{thmnonum}{Theorem}

\newtheorem{thm}{Theorem}
\newtheorem{lem}[thm]{Lemma}

\newcommand{\Z}{\mathbb{Z}}
\newcommand{\R}{\mathbb{R}}
\begin{document}

\title{The number of representations of $n$ as a growing number of squares}
\author{John Holley-Reid}
\email{jholleyreid14@gmail.com}
\author{Jeremy Rouse}
\email{rouseja@wfu.edu}
\begin{abstract}
  Let $r_{k}(n)$ denote the number of representations of the integer $n$
  as a sum of $k$ squares. In this paper, we give an asymptotic for
  $r_{k}(n)$ when $n$ grows linearly with $k$. As a special case, we find that
  \[
     r_{n}(n) \sim \frac{B \cdot A^{n}}{\sqrt{n}},
     \]
with $B \approx 0.2821$ and $A \approx 4.133$.
     
\end{abstract}

\maketitle

\section{Introduction and Statement of Results}

The problem of how many ways a positive integer can be written as a
sum of $k$ squares dates back more than 300 years. In 1640, Fermat
stated (in a letter to Mersenne) that a positive integer $n$ can be
written as the sum of two squares if and only if in the prime
factorization of $n$, the exponents on all primes $p \equiv 3
\pmod{4}$ are even. In 1770, Lagrange proved that every positive
integer is the sum of four squares. In 1834, Jacobi strengthened
Lagrange's theorem and gave a formula for $r_{4}(n)$, the number of
ways that $n$ can be written as a sum of four squares. Jacobi's result
states that
\[
r_{4}(n) = \begin{cases}
  8 \sum_{d | n} d & \text{ if } n \text{ is odd }\\
  24 \sum_{\substack{d | n \\ d \text{ odd }}} d & \text{ if } n \text{ is even. }
  \end{cases}
\]
Jacobi derived this result by considering the Jacobi theta function
\[
\theta(z) = \sum_{n=-\infty}^{\infty} q^{n^{2}}, \quad q = e^{2 \pi i z}
\]
and deriving the relation
\[
\theta(z)^{4} = \sum_{n=0}^{\infty} r_{4}(n) q^{n} = 1 + 8 \sum_{n=1}^{\infty}
\frac{q^{n}}{(1-q^{n})^{2}} - 32 \sum_{n=1}^{\infty} \frac{q^{4n}}{(1-q^{4n})^{2}},
\]
which is equivalent to the formula for $r_{4}(n)$ stated above.

Let $r_{k}(n)$ denote the number of representations of $n$ as the sum
of $k$ squares. Jacobi found formulas for $r_{k}(n)$ expressed in
terms of divisor functions for $k \in \{ 2, 4, 6, 8 \}$ but for values
of $k$, other functions (the so-called cusp forms) are needed. In
1907, Glaisher gave such formulas for even $k \leq 18$ (see
\cite{Glaisher}). Shortly thereafter, Mordell \cite{Mordell} and Hardy \cite{Hardy} applied the circle method to determining asymptotics
for $r_{k}(n)$, and decomposed $r_{k}(n) = \rho_{k}(n) + R_{k}(n)$ as
the sum of the ``singular series'' and an error term. The singular
series has size approximately $n^{(k/2)-1}$ (at least if $k > 4$) and the
error term (when $k$ is even) is $O(d(n) n^{(k/4)-1/2})$ by deep work of
Deligne. Here $d(n)$ is the number of divisors of $n$.

In 2012, the second author determined the implied constant in the estimate
$r_{k}(n) = \rho_{k}(n) + O(d(n) n^{(k/4)-1/2})$. The main result of \cite{Rouse}
is the following.
\begin{thmnonum}
  Suppose that $k$ is a multiple of $4$. If either $k/4$ is odd or $n$ is odd, then we have
\[
  |r_{k}(n) - \rho_{k}(n)| \leq \left(2k + \frac{k (-1)^{k/4}}{(2^{k/2} - 1) B_{k/2}}\right) d(n) n^{\frac{k}{4} - \frac{1}{2}}.
\]
\end{thmnonum}
Here $B_{n}$ is the usual Bernoulli number defined by
\[
\frac{x}{e^{x} - 1} = \sum_{n=0}^{\infty} \frac{B_{n} x^{n}}{n!}.
\]

This result gives a strong bound on $r_{k}(n)$ provided that $n$ is sufficiently
large in terms of $k$. In particular, the main term is larger than the
error term above if $n$ is larger than about $\frac{k^{2}}{\sqrt{\pi e}}$. In
light of this, it is natural to consider the problem of finding
an asymptotic for $r_{k}(n)$ when $k$ and $n$ both grow,
but $n$ is much smaller than $k^{2}$. The main result of this
paper is an asymptotic for when $n$ grows linearly with $k$.

\begin{thm}
\label{main}
  Let $a$ be a positive integer and $b$ be any integer. Then there are constants
  $A$ (depending only on $a$) and $B$ (depending on $a$ and $b$) so that
  \[
    r_{n}(an+b) \sim \frac{B \cdot A^{n}}{\sqrt{n}}.
  \]
\end{thm}
Here $f(n) \sim g(n)$ means that $\lim_{n \to \infty} \frac{f(n)}{g(n)} = 1$.

When $a = 1$ and $b = 0$, the proof produces a value for $A \approx
4.132731376$ and $B \approx 0.28209420367$.
Here is a table of values of $r_{n}(n)$ comparied with $\frac{B \cdot A^{n}}{\sqrt{n}}$.

\begin{tabular}{c|cc}
$n$ & $r_{n}(n)$ & $\frac{B \cdot A^{n}}{\sqrt{n}}$\\
\hline
$10$ & $129064$ & $129648.03$\\
$10^{2}$ & $1.184101 \cdot 10^{60}$ & $1.186074 \cdot 10^{60}$\\
$10^{3}$ & $1.539924 \cdot 10^{614}$ & $1.540180 \cdot 10^{614}$\\
$10^{4}$ & $6.639899 \cdot 10^{6159}$ & $6.640010 \cdot 10^{6159}$\\
$10^{5}$ & $4.657350 \cdot 10^{61620}$ & $4.657358 \cdot 10^{61620}$\\
\end{tabular}

Next, we give a summary of the method we use to prove
Theorem~\ref{main}. We can extract the coefficient
$r_{n}(an+b)$ via
\[
r_{n}(an+b) = \int_{-1/2}^{1/2} q^{-an-b} \theta(x+iy)^{n} \, dx.
\]
To derive asymptotics for this integral, we use the saddle point
method.  The value of the integral above does not depend on $y$, and
we choose $y$ so that $q^{-an-b} \theta(x+iy)^{n}$ has a saddle point
when $x = 0$, which is also the place where the absolute value
of the integrand is maximized. 

In Section 2, we review relevant background and prove a Lemma that gives
an asymptotic for integrals of the type given above. In Section 3, we prove
Theorem~\ref{main} by verifying the hypotheses of the lemma.

\section{Background}
\label{back}

It is not hard to show that for $n \in \Z$,
\[
\int_{-1/2}^{1/2} q^{n} \, dx = \begin{cases}
  1 & \text{ if } n = 0\\
  0 & \text{ if } n \ne 0.
  \end{cases}
\]
This provides a convenient way to extract the coefficient of $q^{n}$ from a generating function. In particular, if $\theta(z) = \sum_{n \in \Z} q^{n^{2}}$,
then $\theta^{n}(z) = \sum_{m=0}^{\infty} r_{n}(m) q^{m}$. Assuming
we are able to switch the infinite sum and integral, we obtain
\begin{equation}
\label{rnform}
  r_{n}(an+b) = \int_{-1/2}^{1/2} q^{-an-b} \theta^{n}(z) \, dx.
\end{equation}
Next, we justify interchanging the sum and integral. In fact we will
also need that for
any non-negative integer $k$, if $f_{m}(z) =
r_{n}(m) q^{m}$, then
\[
\frac{d^{k}}{dz^{k}} \theta^{n}(z) = \sum_{m=0}^{\infty} \frac{d^{k}}{dz^{k}}
f_{m}(z).
\]
Because $\frac{d^{k}}{dz^{k}} f_{m}(z)
= r_{n}(m) (2 \pi i m z)^{k} e^{2 \pi i m z}$, $r_{n}(m)$ is bounded by a polynomial in $m$ (say of degree $r$),
and
\[
\sum_{m=1}^{\infty} m^{r+n} e^{-2 \pi m y}
\]
converges absolutely for any $y > 0$,
it follows that the sequence of partial sums
$\sum_{m=0}^{N} \frac{d^{k}}{dz^{k}} f_{m}(z)$ converges uniformly.
Combining this with the result (see for example Theorem 25.2 on page 185 of \cite{Ross})that if a sequence of functions $g_{n}$ converges
uniformly to $g$, then $\lim_{n \to \infty} \int_{a}^{b} g_{n}(x) \, dx
= \int_{a}^{b} g(x) \, dx$. This yields the desired result.
  
The method of steepest descent (also known as the saddle point method)
is a procedure for obtaining asymptotics for integrals of the form $\int
a(x) b(x)^{n} \, dx$ as $n \to \infty$. For an introduction to this
method, see Chapter 4 of \cite{Miller}. In our case, the method is
quite straightforward, because for any fixed $y$,
the function $\theta(x+iy)$ attains its maximum value at $x = 0$
and the value of $\theta(iy)$ is real. The following
lemma gives us the estimate we desire.

\begin{lem}
\label{lem}
  Suppose that $f(z)$ and $g(z)$ are holomorphic functions and $y \in \R$.
  Suppose that $a$ and $b$ are real numbers with $a < 0 < b$.
  Suppose that for $x \in [a,b]$, $|g(x+iy)| \leq |g(iy)|$ with equality if and only if $x = 0$, $g(iy)$ is real and positive, $f(iy) \ne 0$,
  $g'(iy) = 0$, and $g''(iy)$ is a negative real number. Then we have
  \[
  \int_{a}^{b} f(x+iy) g(x+iy)^{n} \, dx \sim f(iy) g(iy)^{n} \sqrt{\frac{2 \pi g(iy)}{-n g''(iy)}}.
  \]
\end{lem}

Results of this type have appeared in the literature many times before
(see for example Section 5.7, pages 87-89 of \cite{deBruijn}). To keep the
paper self-contained, we provide a complete proof.

\begin{proof}
The function $h(z) = \ln(g(z))$ will be holomorphic in a neighborhood
of $z_{0} = iy$ since $g(iy) > 0$. We consider the Taylor expansion
of $h(z)$ in a neighboorhood of $z = z_{0}$,
\[
h(z) = h(z_{0}) + h'(z_{0}) (z-z_{0}) + \frac{h''(z_{0}) (z-z_{0})^{2}}{2!}
+ E(z).
\]
We have $g'(z_{0}) = \frac{g'(iy)}{g(iy)} = 0$. Moreover, Theorem~8 on
page 125 of \cite{Ahlfors} gives the formula
\[
E(z) = \frac{(z-z_{0})^{3}}{2 \pi i} \int_{\Gamma}
\frac{g(w) \, dw}{(w-z_{0})^{3} (w-z)},
\]
where $\Gamma$ is any simple closed curve contained in the region in which
$h(z)$ is holomorphic that contains $z_{0}$ and $z$. We see then that
\[
|E(z)| \leq |z-z_{0}|^{3} \cdot \left(\text{length of } \Gamma\right)
\cdot \max_{w \in \Gamma} \frac{1}{|w-z_{0}|^{3} |w-z|} \cdot \max_{w \in \Gamma} |h(w)|.
\]
If we require that $|z-z_{0}| < \delta$, we may choose $\Gamma$ to be
a circle of radius $2 \delta$ and it follows that there is a constant $C$
(depending on $\delta$) so that $|E(z)| \leq C |z-z_{0}|^{3}$.

We split up the integral into the contribution near $z_{0}$ (say
the interval $I_{1} = [-1/n^{2/5}, 1/n^{2/5}]$), and the contribution $I_{2}$ away
from $z_{0}$.

For the contribution away from $z_{0}$, once $n$ is large enough,
the maximum value of $f(x+iy) g(x+iy)^{n}$ occurs at either $-1/n^{2/5}$
or $1/n^{2/5}$. The contribution away from $z_{0}$ is hence at most
\[
\max \{ |e^{n h(-n^{-2/5} + iy)} f(-1/n^{4})|, |e^{n h(n^{-2/5} + iy)} f(1/n^{4})| \}.
\]
Using the bound on $E(z)$ above, we see that
\[
n h(\pm n^{-2/5} + iy) = n \ln(g(iy)) + \frac{h''(iy) n^{1/5}}{2!} + O\left(\frac{1}{n^{1/5}}\right).
\]
It follows that as $n$ tends to infinity,
\[
\left|\int_{I_{2}} f(x+iy) g(x+iy)^{n} \, dx\right|
\leq C_{2} f(iy) g(iy)^{n} \cdot e^{h''(iy) n^{1/5}/2}
\]
for some constant $C_{2}$.
Since
\[
h''(iy) = \frac{g''(iy) g(iy) - \left[g'(iy)\right]^{2}}{g(iy)^{2}}
= \frac{g''(iy)}{g(iy)} < 0,
\]
as $n \to \infty$ this contribution is exponentially
smaller than the main contribution.

For the contribution close to $z_{0}$, we consider
\[
\int_{-n^{-2/5}}^{n^{-2/5}} f(x+iy) g(x+iy)^{n} \, dx,
\]
make the change of variables $u = \sqrt{n} x \, dx$, and get
\begin{align*}
  & \int_{-n^{1/10}}^{n^{1/10}} f(u/\sqrt{n} + iy) g(u/\sqrt{n} + iy)^{n} \, \frac{du}{\sqrt{n}}\\
  &= \frac{g(iy)^{n}}{\sqrt{n}} \int_{-n^{1/10}}^{n^{1/10}} f(u/\sqrt{n} + iy)
  e^{u^{2} g''(iy)/2g(iy)} \cdot e^{n E(u/\sqrt{n} + iy)} \, du.
\end{align*}
Fix $\epsilon > 0$. Because $f$ is continuous, it is uniformly continuous
on a short interval surrounding $z_{0}$ and so there is an $N_{1}$ for that
for all $n \geq N_{1}$, and all $u \in [-n^{1/10},n^{1/10}]$, $|f(u/\sqrt{n} + iy) - f(iy)| < \epsilon/2$. Also, for $u \in [-n^{1/10},n^{1/10}]$ we have 
\[
|n E(u/\sqrt{n} + iy)|
\leq n C |u/\sqrt{n}|^{3} \leq C |u|^{3} n^{-1/2} \leq
C n^{3/10} n^{-1/2} = C n^{-2/5}.
\]
It follows that there is some $N_{2}$ so that for $n \geq N_{2}$
$1 - \frac{\epsilon}{2} \leq e^{n E(u/\sqrt{n} + iy)} \leq 1 + \frac{\epsilon}{2}$
for all $u \in [-n^{1/10}, n^{1/10}]$. The triangle inequality then shows that
for $n \geq \max \{N_{1}, N_{2} \}$,
\begin{align*}
& \left|\int_{-n^{-2/5}}^{n^{2/5}} f(x+iy) g(x+iy)^{n} \, dx 
- \frac{f(iy) g(iy)^{n}}{\sqrt{n}} \int_{-n^{1/10}}^{n^{1/10}}
e^{u^{2} g''(iy)/(2g(iy))} \, du\right|\\
&< \epsilon \frac{|f(iy)| g(iy)^{n}}{\sqrt{n}}
\int_{-n^{1/10}}^{n^{1/10}} e^{u^{2} g''(iy)/(2g(iy))} \, du.
\end{align*}
The desired result now follows from making the change of variables
$v = u \sqrt{-g''(iy)/g(iy)} \, du$ and
\[
\lim_{n \to \infty} \int_{-n^{1/10}}^{n^{1/10}}
e^{u^{2} g''(iy)/(2g(iy))} \, du
= \frac{\sqrt{2g(iy)}}{\sqrt{-g''(iy)}} \int_{-\infty}^{\infty} e^{-v^{2}} \, dv
= \frac{\sqrt{2\pi g(iy)}}{\sqrt{-g''(iy)}}.
\]
\end{proof}

\section{Proof of the main result}

In this section, we will prove Theorem~\ref{main} by verifying the hypotheses
of the Lemma~\ref{lem}. We will show that there is a positive real number
$y$ that makes $g'(iy) = 0$. We let $f(z) = q^{-b}$ and
$g(z) = q^{-a} \theta(z)$. It is clear that $f(z)$ is holomorphic,
and (as mentioned in Section~\ref{back}), since we can
differentiate $\theta(z) = 1 + 2 \sum_{n=1}^{\infty} q^{n^{2}}$ termwise,
it follows that $\theta(z)$ is holomorphic as well. Thus, $g(z)$ is holomorphic.

Next, we will show that $|g(x+iy)| \leq |g(iy)|$. We have
\begin{align*}
\left|g(x+iy)\right|&= \left|\sum_{n= -\infty}^{\infty} e^{2 \pi ix(n^{2}-a)} e^{2\pi i(iy)(n^{2}-a)}\right|
\leq \sum_{n= -\infty}^{\infty} \left|e^{2 \pi ix(n^{2}-a)} e^{2\pi i(iy)(n^{2}-a)}\right|\\
&= \sum_{n= -\infty}^{\infty} \left|e^{2 \pi ix(n^{2}-a)}\right| \left|e^{2\pi i(iy)(n^{2}-a)}\right|.
\end{align*}
We know $\left|e^{2\pi i x(n^2-a)}\right| = 1$ to be true. Thus, we have
\[
|g(x+iy)| \leq \sum_{n= -\infty}^{\infty} \left|e^{2\pi i(iy)(n^{2}-a)}\right|
= \sum_{n=-\infty}^{\infty} e^{-2 \pi y(n^{2} - a)} = |g(iy)| = g(iy),
\]
as desired.

Now, we will prove that for $-1/2 \leq x \leq 1/2$,
if $|g(x+iy)| = |g(iy)|$, then $x = 0$. Assume that $|g(x+iy)| = |g(iy)|$
and fix integers $n_{1}$ and $n_{2}$. Using the triangle inequality, we have
\begin{align*}
  |g(x+iy)| &= \left|\sum_{n= -\infty}^{\infty} e^{2 \pi ix(n^{2}-a)} e^{2\pi i(iy)(n^{2}-a)}\right|\\
&= \left|e^{2 \pi i x(n_{1}^{2} - a)} e^{-2 \pi y (n_{1}^{2} - a)}
  + e^{2 \pi i x(n_{2}^{2} - a)} e^{-2 \pi y (n_{2}^{2} - a)}
    + \sum_{n \ne n_{1}, n_{2}} e^{2 \pi i x(n_{1}^{2} - a)} e^{-2 \pi y (n^{2} - a)}\right|\\
    &\leq \left|e^{2 \pi i x(n_{1}^{2} - a)} e^{-2 \pi y (n_{1}^{2} - a)}\right|
      + \left|e^{2 \pi i x(n_{2}^{2} - a)} e^{-2 \pi y (n_{2}^{2} - a)}\right|
        + \sum_{n \ne n_{1}, n_{2}} \left|e^{2 \pi i x(n^{2} - a)} e^{-2 \pi y (n^{2} - a)}\right|\\
        &= \sum_{n=-\infty}^{\infty} |e^{-2 \pi y(n^{2}-a)}| = |g(iy)|.
\end{align*}

Because of the assumption that $|g(x+iy)| = |g(iy)|$, the left hand
side and right hand side of the above inequality are equal, and this
forces all intermediate terms to be equal. Letting $z = e^{2 \pi i
  x(n_{1}^{2} - a)} e^{-2 \pi y (n_{1}^{2} - a)}$ and $w = e^{2 \pi i
  x(n_{2}^{2} - a)} e^{-2 \pi y(n_{2}^{2} - a)}$, we have that $|z+w|
= |z| + |w|$, and a straightforward calcuation shows that $|z+w| = |z|
+ |w|$ forces $z$ and $w$ to have the same phase. This implies that for
any two integers $n_{1}$ and $n_{2}$, $e^{2 \pi i x(n_{1}^{2} - a)} =
e^{2 \pi i x(n_{2}^{2} - a)}$.  Setting $n_{1} = 0$ and $n_{2} = 1$
yields $e^{2 \pi i x} = 1$, and this forces $x = 0$.

It is straightforward to see that $f(iy) \neq 0$ since $f(iy) = e^{2 \pi y b} \ne 0$.

Next, we will show that there is a $y$ so that $g'(iy) = 0$. We have that
\[
g'(iy) = \sum_{n= -\infty}^{\infty} e^{-2\pi y(n^{2}-a)}(-2\pi (n^{2}-a)).
\]
We rewrite the right hand side as
\begin{equation}
\label{posneg}
\sum_{\substack{n \\ n^{2} - a < 0}} (-2 \pi (n^{2} - a)) e^{-2 \pi y (n^{2} - a)}
+ \sum_{\substack{n \\ n^{2} - a \geq 0}} (-2 \pi (n^{2} - a)) e^{-2 \pi y (n^{2} - a)}.
\end{equation}

In the first sum, there are finitely many positive terms, all of which
tend to $\infty$ as $y \to \infty$. As $y \to 0$, we obtain $\sum_{n^{2} - a < 0}
(-2 \pi (n^{2} - a)) > 0$.

In the second sum, the terms are negative and decreasing. It is easy to see that
\[
\frac{-2 \pi e^{2 \pi y}}{(e^{2 \pi y} - 1)^{2}} = \sum_{r=1}^{\infty} -2 \pi r e^{-2 \pi y r} < \sum_{\substack{n \\ n^{2} - a > 0}} (-2 \pi (n^{2} - a)) e^{-2 \pi y (n^{2} - a)} < 0.
\]
Thus, the second sum in \eqref{posneg} tends to zero as $y \to \infty$
and tends to $-\infty$ as $y \to 0$ (since choosing $y$ very small can
make the term $(-2 \pi (n^{2} - a)) e^{-2 \pi y (n^{2} - a)}$
arbitrarily close to $-2 \pi (n^{2} - a)$). It follows that $g'(iy)
\to -\infty$ as $y \to 0$, and $g'(iy) \to \infty$ as $y \to 0$. Since
$g'(iy)$ is continuous, by the intermediate value theorem, there is
some positive real number $y$ for which $g'(iy) = 0$.

Lastly, one can easily see that
\[g''(iy) = (-2\pi i(n^{2}-a))^{2} \sum_{n= -\infty}^{\infty} e^{-2\pi y(n^{2}-a)} < 0.\] 

For the special case that $a = 1$ and $b = 0$, we have $f(z) = 1$
and $g(z) = \sum_{n=-\infty}^{\infty} q^{n^{2} - 1}$. We find that
the value of $y$ that makes $g'(iy) = 0$ is
$y \approx 0.07957745473668$, and this leads to
$A = g(iy) \approx 4.133$ and $B = \sqrt{\frac{2 \pi A}{-g''(iy)}} \approx
0.2821$.

\bibliographystyle{plain}
\bibliography{refs}
\end{document}